\documentclass[12pt]{article}
\usepackage{graphicx}
\usepackage{amssymb}
\usepackage{amsmath}
\numberwithin{equation}{section}

\begin{document}
\author{Andrey Tydnyuk}
\date{November 6, 2006}
\textbf{Rational Solution of the KZ equation (example)}

\begin{center} Andrey Tydnyuk \end{center}
E-mail address:andrey.tydniouk@verizon.net\\
735 Crawford Ave.Brooklyn, NY 11223. USA. \\

\begin{center}{Abstract}\end{center}

We investigate the Knizhnik-Zamolodchikov linear differential
system. The coefficients of this system are rational functions. We
prove that the solution of the KZ system is rational when $k$ is
equal to two and $n$ is equal to three. While doing so, we found the
coefficients of  expansion in a neighborhood of a singular point.


\newpage

\begin{center}{Introduction}\end{center}

The Knizhnik-Zamolodchikov differential system has the form
[1],[2]:
\begin{equation}\frac{dW}{dz}=2A(z)W,\end{equation}

where $A(z)$ and $W(z)$ are $3{\times}3$ matrices,
$z_{1}{\ne}z_{2}$. We suppose that $A(z)$ has the form
\begin{equation}A(z)=\frac{P_{1}}{z-z_{1}}+\frac{P_{2}}{z-z_{2}}.\end{equation}

Here:
\begin{equation}P_{1}=\left[\begin{array}{ccc}
  0 & 1 & 0 \\
  1 & 0 & 0 \\
  0 & 0 & 1 \\
\end{array}\right]\end{equation}

\begin{equation}P_{2}=\left[\begin{array}{ccc}
  0 & 0 & 1 \\
  0 & 1 & 0 \\
  1 & 0 & 0 \\
\end{array}\right]\end{equation}

The matrices $P_{1}$ and $P_{2}$ are connected with the matrix
representation of the symmetric group. In this paper, we consider
the case when $S_{3}$. We prove that in this case the solution of
the Knizhnik-Zamolodchikov is rational. We find the coefficients
of the Laurent expansion in the neighborhood of the point $z_{1}$.
We use the method of L. Sakhnovich  [3].

\section{MAIN NOTIONS}

In a neighborhood of $z_{1}$ the matrix function $A(z)$ can be
represented in the form:
\begin{equation}A(z)=\frac{a_{-1}}{z-z_{1}}+a_{0}+a_{1}(z-z_{1})+...\quad,\end{equation}

where
\begin{equation}a_{-1}=P_{1} ,\quad
a_{r}=(-1)^{r}\frac{P_{2}}{(z_{2}-z_{1})^{r+1}},\quad
r{\geq}0.\end{equation}

\textbf{Proposition 1.1} (see[3])(necessary and sufficient
condition) \emph{If the matrix system}
\begin{equation}[(q+1)I_{3}-2a_{-1}]b_{q+1}=2\sum_{j+\ell=q}a_{j}b_{\ell},\quad-2{\leq}q+1{\leq}2\end{equation}
\emph{has a solution $b_{-2}, \quad b_{-1}, \quad b_{0}, \quad
b_{1}, \quad b_{2}$ and $b_{-2}{\ne}0$. Then system (1) has a
solution}
\begin{equation}W(z)=\sum_{p{\geq}-2}b_{p}(z-z_{1})^{p}, \quad
b_{-2}{\ne}0.\end{equation}
System (1.3) can be written in the
following form
\begin{equation}b_{-2}=I_{3}-P_{1}\end{equation}
\begin{equation}b_{-1}=-2(I_{3}+2P_{1})^{-1}a_{0}b_{-2}\end{equation}
\begin{equation}b_{0}=-P_{1}(a_{0}b_{-1}+a_{1}b_{-2})\end{equation}
\begin{equation}b_{1}=2(I_{3}-2P_{1})^{-1}(a_{0}b_{0}+a_{1}b_{-1}+a_{2}b_{-2})\end{equation}
\begin{equation}(I_{3}-P_{1})b_{2}=(a_{0}b_{1}+a_{1}b_{0}+a_{2}b_{-1}+a_{3}b_{-2})\end{equation}
Direct calculations show that:
\begin{equation}b_{-2}=\left[\begin{array}{ccc}
  1 & -1 & 0 \\
  -1 & 1 & 0 \\
  0 & 0 & 0 \\
\end{array}\right],\end{equation}
\begin{equation}b_{-1}=\frac{1}{-9(z_{2}-z_{1})}\left[%
\begin{array}{ccc}
  -12 & 12 & 0 \\
  6 & -6 & 0 \\
  6 & -6 & 0 \\
\end{array}%
\right],\end{equation}
\begin{equation}b_{0}=\frac{1}{-9(z_{2}-z_{1})^{2}}\left[%
\begin{array}{ccc}
   3 & -3 & 0 \\
  -6 &  6 & 0 \\
   3 & -3 & 0 \\
\end{array}%
\right],\end{equation}
\begin{equation}b_{1}=\frac{1}{-9(z_{2}-z_{1})^{3}}\left[%
\begin{array}{ccc}
  6 & -6 & 0 \\
  6 & -6 & 0 \\
  -12 & 12 & 0 \\
\end{array}%
\right].\end{equation}

It follows from system (1.3) and the relations (1.10)-(1.13) that
\begin{equation}(I_{3}-P_{1})b_{2}=\frac{1}{-9(z_{2}-z_{1})^{4}}\left[%
\begin{array}{ccc}
  1 & -1 & 0 \\
  -1 & 1 & 0 \\
  0 & 0 & 0 \\
\end{array}%
\right].\end{equation} It is easy to see that equation (1.14) has
a solution. Using proposition 1 we obtain the statement.\\

\textbf{Proposition 1.2} \emph{Differential system $(0.1)$ has a
rational fundamental solution.}
\begin{center}{References}\end{center}
1. Chervov A., Talalaev D., KZ equation, G-opers, quantum
Drinfeld-Sokolov reduction and quantum Cayley-Hamilton identity,
arXiv:hep-th/0607250, 2006.\\
2. Etingof P.I., Frenkel I.B.,
Kirillov A.A. jr.,Lectures on Representation Theory and
Knizhnik-Zamolodchikov Equations, Amer.
Math. Society, 1998.\\
3. Sakhnovich L.A., Meromorphic Solutions of Linear Differential
Systems, Painleve Type Functions. Preprint (to appear).
\end{document}